\documentclass[12pt]{article}
\usepackage{amsmath,amsfonts,amssymb,amsthm}
\usepackage{mathtext}
\usepackage{epsfig,enumerate}

\newtheorem*{thm}{Theorem}

\newtheorem{lem}{Lemma}

\begin{document}

\centerline{\bf P.Grinevich, S.Novikov\footnote {P.G.Grinevich,
 Landau Institute for Theoretical Physics, Moscow, Russia\\
e-mail pgg@landau.ac.ru,\\ S.P.Novikov, University of Maryland
at College Park, USA, and Steklov Math Institute, Moscow,
 Russia\\ e-mail novikov@umd.edu}}

\vspace{0.5cm}

\centerline{\Large Spectral Meromorphic Operators and Nonlinear Systems}

\vspace{1cm}

Let us consider ordinary differential linear operators $L=\partial_x^n+\sum_{n\geq i\geq 2}a_i\partial_x^{n-i}$
with $x$-meromorphic coefficients $a_i$ acting in the space of functions of the variable $x\in R\subset C$. We call operator {\bf Spectral Meromorphic or s-meromorphic} at the pole $x_j$ if all solutions $L\psi=\alpha\psi$ are meromorphic at the point $x_j$ for all $\alpha\in C$. We are going to study operators s-meromorphic at all poles $x_j$.

Following functional spaces will be important for our goals. Consider collection of linear subspaces $R_j$ of the space of Laurent polynomials
of the form $\sum_{r_j\leq i< 0}b_iy^i$. Introduce the space of functions $\psi(x)\in F=F_{(\{x_j;R_j\})}$ which are $C^{\infty}$ outside of the poles $x_j$ and meromorphic in the neighborhood of the line $R\subset C$ nearby of the  poles. Its poles are  located
at the fixed discrete set $x_j$. Their Laurent decomposition should have negative parts belonging to the subspaces $R_j$ in the variables
$y=x-x_j$. Positive Laurent coefficients should  satisfy to the following equations: {\bf Every product of functions
from the space $F$ has no terms proportional to $y^{-1}$.}

\begin{lem}For every pair of finitely supported functions $\psi,\phi\in F$ their inner product is well-defined by the formula $$<\psi,\phi>=
\int_R \psi(s)\bar{\phi}(\bar{x})dx$$
The integral here is taken along the real line outside of poles and along any small path surrounding singularity nearby of pole . This inner product is indefinite. Every pole $x_j$ contributes exactly the number of negative squares equal to the dimension of the space of negative Laurent parts $R_j$.
\end{lem}

Proof of lemma easily follows from the property that every product of our functions does not have $y^{-1}$ terms at all poles.\\
As we know, Riemannian surface $\Gamma$ with the marked point $P$, local parameter $z=1/\lambda$ near $P,z(P)=0$, divisor $D$ with degree equal to the genus $g$ and meromorphic function $\alpha(z)$ with the order $k$ pole at the point $P$ only, according to the Burchnall-Chaundy-Krichever (BChK) construction (\cite{K}), define an algebraic (algebro-geometrical)  OD operator $L$ of the order $k$ entering nontrivial rank one commutative ring of the OD operators. This operator is obviously $s$-meromorphic. It is formally self-adjoint (symmetric) if the following $\bf SA$
requirement is satisfied: {\it An antiinvolution is given $\tau: \Gamma\rightarrow\Gamma$ such that $\tau(P)=P,\tau\lambda=-\bar{\lambda}$ and $\tau(D)+D\sim K+2P$. Here $K$ is a canonical divisor of differential forms.}\\
The operator $L$ (if it is singular) defines a space of functions $F^L=F$ by its singularities (see above).For periodic operators we define  subspaces
$F^L_{\kappa}\subset F^L$ for $|\kappa|=1$ such that $f(x+T)=\kappa f(x)$ for $f\in F^L_{\kappa}$. In the spaces $F^L_{\kappa}$ the number of negative squares of inner product is
finite. It does not depend on $\kappa$ and equal to $\sum_{a\leq x_j\leq a+T} \dim R_j$. Every singular point $x_j$ contributes exactly $dim R_j$
negative squares. \\
By definition, the quasimomentum differential $dp$ is meromorphic and such that $dp=d\lambda +Reg$ near $P$. All its periods  $\oint dp$ along every cycle are real. We choose  imaginary part $p_I$
as a one-valued real function such that $\tau p_I=-p_I$ because $\tau dp=\bar{dp}$. The level $p_I=0$ we call a {\bf Canonical Contour}.
\begin{thm}
I.In the spaces $F^L$ for the algebraic operators $L$ satisfying to the SA requirement, the inner product defined above is well-defined and indefinite. Our operator is self-adjoint with respect to this product  in the spaces $F^L$ .\\
II. Assume now that  the Baker-Akhiezer eigenfunctions $\psi(z,x)$ corresponding to the contour $z\in C_0$ give complete family in the space $F^L$.
  The spectrum coincides with projection $\alpha(C_0)$ of the canonical contour to the complex plane. The space of functions with indefinite inner product can be realized by the smooth functions at the contour $C_0$ using spectral measure $\Omega$--see below. The number of negative squares in the spaces $F^L_{\kappa}$ is equal to the number of negative points for the spectral measure $\Omega|_{z_j}<0 $ where $e^{ip(z_j)T}=\kappa$.
  In particular, this number is an integral of time dynamics generated by the Gelfand-Dikii systems with Lax representation $\dot L= [L,A]$ such that the orders of operators are relatively prime and $L$ is (hermitian) symmetric.
\end{thm}
To prove this theorem, we consider the Baker-Akhiezer function $\psi(\lambda,x)\sim e^{\lambda x}$ near $P$, with the divisor of poles $D=(\gamma_1,...,\gamma_g)$ and dual 1-form $\psi^*(\mu,x)d\mu\sim e^{-\mu x}d\mu$ near $P$ with zeroes $D=(\gamma_1,...,\gamma_g)$. In terms of the dual scalar Baker-Akhiezer function $\psi^+$ we have $\psi^*(\mu,x)d\mu)=\psi^+(\mu,x)\Omega$ where $\Omega=(1+o(1))d\mu$ near $P$  with zeroes
$D+\tau D\sim K+2P$, $K$ is the canonical divisor of differential forms. $\Omega$ serves as a  {\bf spectral measure} at the contour $C_0$.
We have $L^*\psi^*=\mu\psi^*$ where $L^*$ is an operator formally adjoint to $L$ (applying complex conjugation, we obtain similar result for the hermitian adjoint operators). This fact never was formally proved in the literature but known at the ''oral level''  among experts. The authors of the work \cite{GO}
studied the Cauchy kernel
on the algebraic curve $\Gamma$ such that $\omega(\lambda,\mu,x)\sim d\mu/(\mu-\lambda)+O(1)$ near the diagonal $\mu=\lambda$, $\omega\sim e^{\lambda x}o(1)$ for $\lambda\rightarrow P$, $\omega\sim e^{-\mu x}o(1)d\mu$ for $\mu\rightarrow P$.  The equality was found $\partial_x\omega=-\psi(\lambda,x)\psi^*(\mu,x)d\mu$. This fact is crucial here. It follows that every product of eigenfunctions have residue term
(in the variable $x$) equal to zero. Therefore our inner product (above) is well-defined. For the given unimodular Bloch multiplier $\kappa$
the number of negative squares of inner product is an invariant of Riemann surface $\Gamma$--for the space of smooth functions at the contour $C_0$
with spectral measure $\Omega$: The equation $e^{ip(z)T}=\kappa$  gives a countable discrete set of points $z_l\in C_0$. The number of point with negative sign $\Omega|_{z_l}$ is finite and equal to the number of negative squares of the inner product at the contour $C_0$. This number remains unchanged with time and $\kappa$. To prove this fact we  consider  $\kappa$ dependence of the equation $e^{ip(z)T}=\kappa$. Something might happen with sign of $\Omega|_{C_0}$ if our divisor $\tau D(t)+D(t)$ crosses also $C_0$ at the same point. Let $\tau(z_l)=z_l$. In this case we have $\tau D$ and $D$ coincide at this point. Therefore the zero of $\Omega$ has even order> and sign of $\Omega$ remains unchanged The second possibility is that $\tau z_l\neq z_l$. This case is trivial: Inner product is and was indefinite here with equal number of  positive and  negative squares before and after crossing. All set near $P$ at the contour $C_0$ always give positive inner product. So inner product at the contour $C_0$ does not depend on $\kappa$
(if all solutions $z_l$ has definite sign). Therefore it does not depend on time. Assuming that our eigenfunctions are complete in the space $F^L_{\kappa}$, we obtain the proof of our theorem.\\
{\bf Example}(see \cite{GN}). Let $n=2$ and $L=-\partial_x^2+u$. For $s$-meromorphic potentials singularities should have form $u=r_j(r_j+1)/(x-x_j)^2+Regular$
where regular part has a form $Regular=\sum_{i\leq r_j-1}c_{ij} (x-x_j)^{2i}+O((x-x_j)^{2r_j})$. The space $F^L$ consists of functions with Laurent parts
$f_(x)=b_1y^{-r_j}+b_2y^{-r_j+2}+....+b_{r_j}y^{r_j}+O(y^{r_j+1})$  where $y=x-x_j$. Every pole contributes exactly $[(r_j+1)/2]$ negative squares to the spaces like
$F^L_{\kappa}$ or rapidly decreasing functions on the $x$-line. For genus one and rombic lattice one component of the spectrum $\alpha(C_0)$ is complex.
This case never appeared in the Hermit-Lame' spectral theory--only rectangular lattice appeared in the classical theory of that operator.\\
Following questions 1-2 will be answered in the next work:\\
1.For $s$-meromorphic operators the adjoint operators are also $s$-meromorphic.\\
2.The Baker-Akhiezer functions give a complete family in the spaces $F^L,F^L_{\kappa}$\\
{\bf Problem}. Classify all singularities of $s$-meromorphic (in particular, algebraic rank one) operators. Prove that these two classes of singularities coincide. In order to find the spaces $R_j$ of Laurent polynomials we need to consider the Baker-Akhiezer function of the whole KP hierarchy $\psi(x,y,t,...)$ and find its $x$-singularities depending on other variables as from parameters.\\
For the second order operators $n=2$ all these problems are solved (see \cite{GN}). {\bf In process of  real time evolution singularity may jump
preserving the number of negative squares of the inner product}. For example, for  the singular 2-soliton KdV solution with one singularity of the type $2/(x-x_0(t))^2$
 might jump to the type $6/(x-x_0)^2$
at some isolated moments of time.

{\bf Nonstationary (1+1) Periodic Schrodinger Operator}. This theory works also for the Nonstationary Schrodinger Operator $L=i\partial_y+\partial_x^2+ u(x,y)$
acting in the space of functions $F(x,y)$ as well as in its ''physical subspace''  $F_0=\{Lf=0\}$. Same arguments as above, using results of \cite{GO}, lead to conclusion that
operator $L$ is symmetric in respect to the same indefinite inner product $<f,g>=\int f(x,y)\bar{g}(\bar{x},y)dxdy$ if singularities of $u$ are the same as for algebrogeometrical (Krichever)
solutions to the KP-I equation. We hope to prove soon the completeness of this basis. In the regular case the completeness was proved by Krichever \cite{K2} for periodic solutions . In the work \cite{DN}
 regularity criteria for the real solutions to KP-I is in fact associated with positivity of inner product of corresponding Baker-Akhiezer functions with spectral measure
 along the canonical contour  on the Riemann surface. Completeness (if it is true) implies that these inner product is equivalent to our.

{\bf Abstract:}

We study here class of 1D {\bf spectral-meromorphic (s-meromorphic)}
OD  operators $L=\partial_x^n+\sum_{n-2\geq i\geq 0}a_{n-2-i}\partial_x^i$
with meromorphic coefficients $a_j$ near $x\in R$
such that all eigenfunctions $L\psi=\alpha\psi$ are $x$--meromorphic near $x\in R$ for all
$\alpha$.
Symmetric  $s$-meromorphic operators are
self-adjoint with respect to indefinite inner product well-defined for some special spaces
of singular functions.
In particular, all algebraic operators $L$--i.e. operators entering
Burchnall-Chaundy-Krichever (BChK)
{\bf rank one} commutative rings --are s-meromorphic.
For $KdV$ system corresponding algebraic operator  $L=-\partial_x^2+u(x,t)$ is called
singular finite gap, singular soliton or algebrogeometric Schrodinger operator. This
special case was already studied by the present authors  in the recent works (see
\cite{GN} and references therein). New results are added in this version about (1+1)
Nonstationary Schrodinger operators

\end{document}